\documentclass[12pt]{amsart}   
\usepackage{amsmath,amsthm,amssymb,marvosym}
\usepackage[english]{babel}
\usepackage[autostyle]{csquotes}

\usepackage{graphicx,subfig,wrapfig}
\usepackage{textcase}
\usepackage{enumitem, ragged2e}
\usepackage{lipsum,hyperref}
\newtheorem{theorem}{Theorem}[section]
\newtheorem{Lemma}[theorem]{Lemma}

\newtheorem{Definition}[theorem]{Definition}

\newtheorem{Remark}[theorem]{Remark}
\newtheorem{Example}[theorem]{Example}
\newdimen\slantmathcorr
\def\oversl#1{
\setbox0=\hbox{$#1$}
\slantmathcorr=\wd0
\hskip 0.2\slantmathcorr \overline{\hbox to 0.8\wd0{
\vphantom{\hbox{$#1$}}}}
\hskip-\wd0\hbox{$#1$}
}

\newtheoremstyle{Proof}{}{}{}{}{\bfseries}{.}{.5em}{}
\theoremstyle{Proof}
\newtheorem*{Proof}{Proof}

\usepackage{ifthen}

\makeatletter
\renewcommand\@cite[2]{%
~[#1]\ifthenelse{\boolean{@tempswa}}
{(\nolinebreak[3]#2)}{}}
\renewcommand\@biblabel[1]{#1.}
\makeatother
\begin{document} 
\title[\tiny{T\lowercase{opologically stable and} $\beta$-\lowercase{persistent points of group actions}}]{Topologically stable and $\beta$-persistent points of group actions}
\author[Abdul Gaffar Khan \lowercase{and} Tarun Das]{Abdul Gaffar Khan$^{1,2}$ \lowercase{and} Tarun Das$^{1,3}$ \\ \tiny{D\lowercase{epartment of} M\lowercase{athematics}, U\lowercase{niversity of} D\lowercase{elhi}, D\lowercase{elhi}, I\lowercase{ndia.}$^{1}$}\\
\lowercase{gaffarkhan18@gmail.com}$^{2}$ \lowercase{and tarukd@gmail.com}$^{3}$}
%\keywords{Topological stability, Persistent, Borel measure.\vspace*{0.08cm}\\ 
%\hspace*{0.27cm} \textit{2020 Mathematics Subject Classification.} Primary: 37C85; Secondary: 37B25.\\
%\vspace*{0.01cm}
%\Letter{Tarun Das} \\
%\vspace*{0.08cm}
%tarukd@gmail.com \\
%\vspace*{0.01cm}
%Abdul Gaffar Khan \\
%\vspace*{0.08cm}
%gaffarkhan18@gmail.com\\
%$^{1}$\textit{Department of Mathematics, Faculty of Mathematical Sciences,} 
%\textit{University of Delhi, Delhi, India.} 
%}

\begin{abstract}
In this paper, we introduce topologically stable points, $\beta$-persistent points, $\beta$-persistent property, $\beta$-persistent measures and almost $\beta$-persistent measures for first countable Hausdorff group actions of compact metric spaces. We prove that the set of all $\beta$-persistent points is measurable and it is closed if the action is equicontinuous. We also prove that the set of all $\beta$-persistent measures is a convex set and every almost $\beta$-persistent measure is a $\beta$-persistent measure. Finally, we prove that every equicontinuous pointwise topologically stable first countable Hausdorff group action of a compact metric space is $\beta$-persistent. In particular, every equicontinuous pointwise topologically stable flow is $\beta$-persistent. 
\end{abstract} 
\maketitle 
\vspace*{-0.5cm}
\begin{flushleft}
\footnotesize{\textit{Key words and phrases.} Topological stability, Persistent, Borel measure.}\\
\end{flushleft}

\section{Introduction}
In \cite{WO}, author has studied topologically stable homeomorphisms of compact metric spaces. In \cite{LP}, Lewowicz called a homeomorphism to be persistent if every orbit of the map can be seen through some actual orbit of every small enough perturbed system. In the current literature, persistent property is also referred as $\beta$-persistent property. Every topologically stable homeomorphism of a compact manifold is persistent but this need not be true when the phase space is a compact metric space \cite[Example 1]{SKO}. Recently, in  \cite{KLMP} authors have introduced pointwise topologically stable homeomorphisms. First motivation of this paper comes from the relationship obtained between pointwise topological stability and persistent property of a homeomorphism in \cite{DLMP, KDP}. Precisely, authors have proved that every equicontinuous pointwise topologically stable homeomorphism of a compact metric space is persistent. 
\par 
\vspace*{0.2cm}

In \cite{CLT}, authors have introduced topologically stable finitely generated group actions of compact metric spaces and in \cite{DJLP}, authors have studied pointwise topologically stable finitely generated group actions. 
In \cite{CG}, author has introduced $GH$-stable countable group actions which extends the notion of topologically $GH$-stable homeomorphisms \cite{ART} and $GH$-stable finitely generated group actions \cite{DLMG, KDDG} in more general setting where $GH$ stands for Gromov-Hausdorff distance. 
Second motivation of this paper comes from the study of such stability for countable group actions and from the notion of $\beta$-persistent property of finitely generated group actions \cite{ALLP}. 
We combine the idea of distance used in \cite{CG} with notions studied in \cite{ALLP} to introduce and study pointwise topological stability and persistent property for first countable Hausdorff group actions where we assume that the group admits a dense countable subgroup.
First we introduce topologically stable points and $\beta$-persistent property for first countable Hausdorff group actions. Then we extend the relationship between pointwise topologically stable homeomorphisms and its $\beta$-persistent property obtained in \cite{DLMP, KDP} to such actions. Precisely, we prove the following. 

\begin{theorem}\label{T1.1}
Let $\Phi: G\times X\rightarrow X$ be a continuous first countable Hausdorff group action of a compact metric space $(X, d)$, where $G$ admits a dense countable subgroup. If $\Phi$ is equicontinuous pointwise topologically stable, then $\Phi$ is $\beta$-persistent.
\end{theorem}

\begin{Remark}\label{C1.2}
Note that the union of the set of all $n^{th}$-roots of unity, where $n\in \mathbb{N}$, is a dense countable subgroup in the unit circle $\mathbb{S}^{1}$ and the set of all rational numbers is a dense countable subgroup in $\mathbb{R}$. 
Therefore, Theorem \ref{T1.1} can be applied to $G$-actions, where $G = \mathbb{R}$ or $\mathbb{S}^{1}$. To study topologically stable flows, please refer to \cite{HM}. 
\end{Remark}

This paper is distributed as follows. Section 2 gives necessary preliminaries required for the remaining section. 
In Section 3, we introduce topologically stable points, $\beta$-persistent points, $\beta$-persistent property and $\beta$-persistent measures for first countable Hausdorff group actions. Then we prove Theorem \ref{T1.1}.

\section{Preliminaries}
Throughout this paper, $G$ denotes a first countable Hausdorff group, $H$ denotes a countable dense subgroup of $G$ and $(X, d)$ denotes a compact metric space. The set of all natural numbers is denoted by $\mathbb{N}$. For a given $\delta > 0$ and for each $x\in X$, we denote $B(x, \delta)= \lbrace y\in X \mid d(x, y) < \delta \rbrace$ and $B[x, \delta]= \lbrace y\in X \mid d(x, y) \leq \delta \rbrace$.
\par
\vspace*{0.2cm}

We say that $G$ is a finitely generated group if there exists a finite symmetric generating set $S=\lbrace s_{i} \mid 1\leq i\leq n\rbrace$ of $G$, for some $n\in \mathbb{N}$ i.e. $S$ is a finite subset of $G$, $s\in S$ if and only if $s^{-1}\in S$ and every element of $G$ can be expressed as a combination of finitely many elements of $S$.  
\par
\vspace*{0.2cm}

A group action of $X$ with respect to the group $G$ is a continuous map $\Phi : G\times X\rightarrow X$ such that:\\
(i) For each $g\in G$, $\Phi_{g}(.) = \Phi(g, .)$ is a homeomorphism of $X$.\\
(ii) For the identity element $e\in G$, $\Phi_{e}(x) = x$, for each $x\in X$.\\
(iii) $\Phi_{g}(\Phi_{h}(x)) = \Phi_{gh}(x)$, for each pair $g, h\in G$ and for each $x\in X$.
\par
\vspace*{0.2cm}

The class of all group actions of $X$ with respect to $G$ is denoted by $\operatorname{Act(G,X)}$. The orbit of a point $x\in X$ under $\Phi\in \operatorname{Act(G, X)}$ is given by $\mathcal{O}_{\Phi}(x) = \lbrace \Phi_{g}(x) \mid g\in G\rbrace$. For a subset $J$ of $G$, we define the distance $d_{J}$ on $\operatorname{Act(G, X)}$ by $d_{J}(\Phi, \Psi) = \sup\limits_{j\in J, x\in X} d(\Phi_{j}(x), \Psi_{j}(x))$, for each pair $\Phi, \Psi\in \operatorname{Act(G, X)}$. Note that $d_{K}(\Phi, \Psi) \leq d_{G}(\Phi, \Psi)$, for each pair $\Phi, \Psi\in \operatorname{Act(G, X)}$ and for every subset $K$ of $G$. 

\begin{Remark}\label{R2.1}
Let $\Phi, \Psi\in \operatorname{Act(G, X)}$ and $H_{1}, H_{2}$ be two dense countable subgroups of $G$. Since every first countable Hausdorff group is metrizable, we can use continuity of $\Phi$ and $\Psi$ to conclude that $d_{H_{1}}(\Phi, \Psi) = d_{H_{2}}(\Phi, \Psi)$.
\end{Remark}

Let $(Y, p)$ be a metric space and $H: X\rightarrow Y$ be a homeomorphism. For each $\Phi \in \operatorname{Act(G, X)}$, we define the action $H \Phi H^{-1}\in \operatorname{Act(G, Y)}$ by $(H \Phi H^{-1})_{g} = H\circ \Phi_{g}\circ H^{-1}$, for each $g\in G$. 
\par
\vspace*{0.2cm}

Let $\Phi \in \operatorname{Act(G, X)}$ and $H$ be a dense countable subgroup of $G$. We say that $\Phi$ is equicontinuous (with respect to $H$) if for each $\epsilon > 0$, there exists a $\delta > 0$ such that if $x,y\in X$ satisfy $d(x, y) \leq \delta$, then $d(\Phi_{h}(x), \Phi_{h}(y)) \leq \epsilon$, for each $h\in H$. 
We say that a point $z\in X$ is an equicontinuous point of $\Phi$ (with respect to $H$) if for each $\epsilon > 0$, there exists a $\delta > 0$ depending on $\epsilon$ and $z$ such that if $y\in X$ satisfies $d(z, y) \leq \delta$, then $d(\Phi_{h}(z), \Phi_{h}(y)) \leq \epsilon$, for each $h\in H$. The set of all equicontinuous points of $\Phi$ is denoted by $E(\Phi)$.
%\par
%\vspace*{0.2cm}

\begin{Remark}\label{R2.2}
Let $\Phi \in \operatorname{Act(G, X)}$ and $H$ be a dense countable subgroup of $G$. Since $G$ is a first countable Hausdorff group and $\Phi$ is continuous, we can get that $\Phi$ is equicontinuous (with respect to $H$) if and only if for each $\epsilon > 0$, there exists a $\delta > 0$ such that if $x,y\in X$ satisfy $d(x, y) \leq \delta$, then $d(\Phi_{g}(x), \Phi_{g}(y)) \leq \epsilon$, for each $g\in G$. Also,  $z\in X$ is an equicontinuous point of $\Phi$ (with respect to $H$) if and only if for each $\epsilon > 0$, there exists a $\delta > 0$ depending on $\epsilon$ and $z$ such that if $y\in X$ satisfies $d(x, y) \leq \delta$, then $d(\Phi_{g}(z), \Phi_{g}(y)) \leq \epsilon$, for each $g\in G$.
\end{Remark}

Let $f:X\rightarrow X$ be a homeomorphism. We say that a point $x\in X$ is a topologically stable point of $f$ if for each $\epsilon > 0$, there exists a $\delta > 0$ such that for each homeomorphism $g:X\rightarrow X$ satisfying $\sup\limits_{y\in X} d(f(y), g(y))\leq \delta$, there exists a continuous map $h : \overline{\mathcal{O}_{g}(x)} \rightarrow X$ such that $f\circ h = h\circ g$ and $d(h(z), z) \leq \epsilon$, for each $z\in \overline{\mathcal{O}_{g}(x)}$ \cite{KLMP}. 
\par
\vspace*{0.2cm}

Let $\Phi\in \operatorname{Act(G,X)}$ and $S$ be a finite symmetric generating set of $G$. We say that $\Phi$ is topologically stable (with respect to $S$) if for each $\epsilon >0$, there exists a $\delta > 0$ such that for each $\Psi\in \operatorname{Act(G, X)}$ satisfying $d_{S}(\Phi, \Psi) \leq \delta$, there exists a continuous map $h:X\rightarrow X$ such that $\Phi_{g}\circ h = h\circ \Psi_{g}$, for each $g\in G$ and $d(h(x), x) \leq \epsilon$, for each $x\in X$. Recall that the topological stability of $\Phi$ is independent of the choice of a generator of $G$. 
We say that $\Phi$ is topologically stable if $\Phi$ is topologically stable with respect to some symmetric generating set of $G$  \cite{CLT}. 
We say that a point $x\in X$ is a topologically stable point of $\Phi$ (with respect to $S$) if for each $\epsilon > 0$, there exists a $\delta > 0$ such that for each $\Psi\in \operatorname{Act(G, X)}$ satisfying $d_{S}(\Phi, \Psi) \leq \delta$, there exists a continuous map $h : \overline{\mathcal{O}_{\Psi}(x)} \rightarrow X$ such that $\Phi_{g}\circ h = h\circ \Psi_{g}$, for each $g\in G$ and $d(h(z), z) \leq \epsilon$, for each $z\in \overline{\mathcal{O}_{\Psi}(x)}$.  
We say that $x$ is a topologically stable point of $\Phi$ if $x$ is a topologically stable point of $\Phi$ with respect to some symmetric generating set of $G$. The set of all topologically stable points of $\Phi$ is denoted by $T_{G}(\Phi)$. We say that $\Phi$ is pointwise topologically stable if $T_{G}(\Phi) = X$. Recall that if $\Phi$ is topologically stable, then $\Phi$ is pointwise topologically stable \cite{DJLP}.
%\par
%\vspace*{0.2cm}

Let $\Phi\in \operatorname{Act(G, X)}$ and $S$ be a finite symmetric generating set of $G$. We say that $\Phi$ is $\beta$-persistent (with respect to $S$) if for each $\epsilon > 0$, there exists a $\delta > 0$ such that for each $\Psi\in \operatorname{Act(G, X)}$ satisfying $d_{S}(\Phi, \Psi) \leq \delta$ and for each $x\in X$, there exists a $y\in X$ such that $d(\Phi_{g}(x), \Psi_{g}(y)) \leq \epsilon$, for each $g\in G$. Recall that the $\beta$-persistent property of $\Phi$ is independent of the choice of a generator of $G$. We say that $\Phi$ is $\beta$-persistent if $\Phi$ is $\beta$-persistent with respect to some symmetric generating set of $G$ \cite{ALLP}.
\par
\vspace*{0.2cm}

A Borel probability measure $\mu$ of $X$ is a $\sigma$-additive measure defined on the Borel $\sigma$-algebra $\mathcal{B}(X)$ of $X$  such that $\mu(X)  = 1$. The set of all Borel probability measures of $X$ is denoted by $\mathcal{M}(X)$ and it is assumed to be equipped with the weak$^{*}$ topology. We say that $\mu$ is supported on a subset $B$ of $X$ if the support of the measure $\mu$, denoted by $supp(\mu)$, satisfies $supp(\mu) \subseteq B$. We denote the Dirac measure supported on $\lbrace x\rbrace \subseteq X$ by $m_x$ i.e. $m_x(A) = 0$, if $x\notin A$ and $m_{x}(A) = 1$, if $x\in A$, for every subset $A$ of $X$.

\section{Proof of Theorem \ref{T1.1}}
In this section, we first define topologically stable points for first countable Hausdorff group actions and compare it with topologically stable points of a finitely generated group action  \cite{DJLP}. 
Then we define $\beta$-persistent points and $\beta$-persistent measures for such actions to prove Theorem \ref{T1.1}. 

\begin{Definition}\label{D3.1}
Let $\Phi\in \operatorname{Act(G,X)}$ and $H$ be a dense countable subgroup of $G$. We say that a point $x\in X$ is a topologically stable point of $\Phi$ (with respect to $H$) if for each $\epsilon >0$, there exists a $\delta > 0$ such that for each $\Psi\in \operatorname{Act(G, X)}$ satisfying $d_{H}(\Phi, \Psi) \leq \delta$, there exists a continuous map $h: \overline{\mathcal{O}_{\Psi}(x)}\rightarrow X$ such that $\Phi_{g}\circ h = h\circ \Psi_{g}$, for each $g\in G$ and $d(h(y), y) \leq \epsilon$, for each $y\in \overline{\mathcal{O}_{\Psi}(x)}$. The set of all topologically stable points of $\Phi$ is denoted by $T(\Phi)$. We say that $\Phi$ is pointwise topologically stable if $T(\Phi) = X$.
\end{Definition}

\begin{Remark}\label{E1R3.2}
Let $\Phi\in \operatorname{Act(G, X)}$ and $H_{1}, H_{2}$ be dense countable subgroups of $G$. From Remark \ref{R2.1}, we get that $x$ is a topologically stable point of $\Phi$ (with respect to $H_{1}$) if and only if $x$ is a topologically stable point of $\Phi$ (with respect to $H_{2}$). Thus the notion of topologically stable points is independent of the choice of dense countable subgroups. 
\end{Remark}

\begin{Remark}\label{R3.2}
Let $G$ be a finitely generated group with finite symmetric generating set $S$ and $\Phi\in \operatorname{Act(G, X)}$. Since $d_{S}(\Phi, \Psi) \leq d_{G}(\Phi, \Psi)$, for each $\Psi\in \operatorname{Act(G, X)}$, we use Remark \ref{E1R3.2} to get that $T_{G}(\Phi)\subseteq T(\Phi)$. Therefore if $\Phi$ is topologically stable, then $T_{G}(\Phi) = X = T(\Phi)$ and hence $\Phi$ is pointwise topologically stable in the sense of Definition \ref{D3.1}. 
\end{Remark}

\begin{theorem}
Let $\Phi \in \operatorname{Act(G, X)}$. Then $T(R \Phi R^{-1}) = R(T(\Phi))$, for each homeomorphism $R:X \rightarrow Y$.
\label{T3.3}
\end{theorem}
\begin{Proof}
Choose an $x\in T(\Phi)$. For a given $\epsilon > 0$, choose an $\eta > 0$ by uniform continuity of $R$. For this $\eta$, choose a $\delta > 0$ by the definition of topologically stable point $x$ of $\Phi$. For this $\delta$, choose a $\gamma > 0$ by uniform continuity of $R^{-1}$. Choose a $\Psi \in \operatorname{Act(G, Y)}$ satisfying $p_{H}(R \Phi R^{-1}, \Psi) = \sup\limits_{h\in H, x\in X}p(R \Phi_{h} R^{-1}(x), \Psi_{h}(x)) \leq \gamma$. 
Clearly $d_{H}(\Phi, R^{-1} \Psi R) \leq \delta$ and hence there exists a continuous map $h : \overline{\mathcal{O}_{R^{-1} \Psi R}(x)} \rightarrow X$ such that $\Phi_{g}\circ h = h\circ R^{-1}\Psi_{g} R$, for each $g\in G$ and $d(h(y), y) \leq \epsilon$, for each $y\in \overline{\mathcal{O}_{R^{-1}\Psi R}(x)}$ . 
Note that the map $h' : \overline{\mathcal{O}_{\Psi}(R(x))}\rightarrow Y$ defined by $h' = R h R^{-1}$ is a well defined continuous map such that $R \Phi_{g} R^{-1} h = h \Psi_{g}$, for each $g\in G$ and $p(h'(z), z) \leq \epsilon$, for each $z \in \overline{\mathcal{O}_{\Psi}(R(x))}$. 
Since $\epsilon$ and $x$ are chosen arbitrarily, we get that $R(T(\Phi)) \subseteq T(R \Phi R^{-1})$. Replace $X$ by $Y$, $Y$ by $X$, $R$ by $R^{-1}$ and $\Phi$ by $R \Phi R^{-1}$ in the last inclusion to complete the proof. \qed
\end{Proof}

\begin{Definition}
Let $\Phi\in \operatorname{Act(G, X)}$ and $H$ be a dense countable subgroup of $G$. We say that $\Phi$ is $\beta$-persistent (with respect to $H$) through a subset $B$ of $X$ if for each $\epsilon > 0$, there exists a $\delta > 0$ such that for each $\Psi\in \operatorname{Act(G, X)}$ satisfying $d_{H}(\Phi, \Psi) \leq \delta$ and for each $x\in B$, there exists a $y\in X$ such that $d(\Phi_{g}(x), \Psi_{g}(y)) \leq \epsilon$, for each $g\in G$. 
We say that $\Phi$ is $\beta$-persistent (with respect to $H$) if $\Phi$ is $\beta$-persistent (with respect to $H$) through $X$.
We say that a point $x\in X$ is a $\beta$-persistent point of $\Phi$ (with respect to $H$) if $\Phi$ is $\beta$-persistent (with respect to $H$) through $\lbrace x\rbrace$. The set of all $\beta$-persistent points of $\Phi$ is denoted by $P(\Phi)$. We say that $\Phi$ is pointwise $\beta$-persistent if $P(\Phi) = X$.
\label{D3.4}
\end{Definition}

\begin{Remark}\label{E1R3.5.1}
Let $\Phi\in \operatorname{Act(G, X)}$ and $H_{1}, H_{2}$ be dense countable subgroups of $G$. From Remark \ref{R2.1}, we get that $\Phi$ is $\beta$-persistent (with respect to $H_{1}$) through a subset $B$ of $X$ if and only if $\Phi$ is $\beta$-persistent (with respect to $H_{2}$) through $B$. In particular, the notion of $\beta$-persistent property of $\Phi$ and the notion  of $\beta$-persistent point of $\Phi$ are independent of the choice of dense countable subgroups. 
\end{Remark}

\begin{Remark}\label{R3.5}
Let $G$ be a finitely generated group with finite symmetric generating set $S$ and $\Phi\in \operatorname{Act(G, X)}$. Since $d_{S}(\Phi, \Psi) \leq d_{G}(\Phi, \Psi)$, for each $\Psi\in \operatorname{Act(G, X)}$, we use Remark \ref{E1R3.5.1} to get that if $\Phi$ is $\beta$-persistent (with respect to $S$), then $\Phi$ is $\beta$-persistent in the sense of Definition \ref{D3.4}.
\end{Remark}

\begin{Remark}\label{R3.6}
Let $\Phi\in \operatorname{Act(G, X)}$ and $R:X\rightarrow Y$ be a homeomorphism. Then the following statements are true:
\begin{enumerate}
\item[(1)] If $\Phi$ is $\beta$-persistent, then $\Phi$ is pointwise $\beta$-persistent.
\item[(2)] $P(R\Phi R^{-1}) = R(P(\Phi))$.
\end{enumerate}
\end{Remark}

Let $\Phi , \Psi \in \operatorname{Act(G, X)}$. Then for each $\epsilon > 0$ and for each $x\in X$, we denote $\Gamma_{\epsilon}^{x}(\Phi , \Psi) = \bigcap\limits_{g\in G} \Psi_{g^{-1}}(B[\Phi_{g}(x), \epsilon]) = \lbrace y\in X \mid d(\Phi_{g}(x), \Psi_{g}(y)) \leq \epsilon,$ for each $g\in G \rbrace$. We define $B(\epsilon , \Phi , \Psi) = \lbrace x\in X \mid \Gamma_{\epsilon}^{x}(\Phi , \Psi) \neq \phi \rbrace$.

\begin{Lemma}
Let $\Phi , \Psi \in \operatorname{Act(G, X)}$. Then $B(\epsilon , \Phi, \Psi)$ is a compact subset of $X$, for each $\epsilon > 0$.
\label{L3.7}
\end{Lemma}
\begin{Proof}
Since $X$ is compact, it is sufficient to show that $B(\epsilon, \Phi, \Psi)$ is a closed subset of $X$, for each $\epsilon > 0$. Fix an $\epsilon > 0$ and choose a sequence $\lbrace x_{i}\rbrace_{i\in \mathbb{N}}$ in $B(\epsilon, \Phi, \Psi)$ such that $x_{i}\rightarrow x$, for some $x\in X$. Then there exists a sequence $\lbrace y_{i}\rbrace_{i\in \mathbb{N}}$ in $X$ such that $d(\Phi_{g}(x_{i}), \Psi_{g}(y_{i})) \leq \epsilon$, for each $g\in G$ and for each $i\in \mathbb{N}$. Since $X$ is a compact metric space, we can assume that $y_{i}\rightarrow y$, for some $y\in X$. Note that $y\in \Gamma_{\epsilon}^{x}(\Phi, \Psi)$ and hence $x\in B(\epsilon, \Phi, \Psi)$ implying that $B(\epsilon, \Phi, \Psi)$ is a closed subset of $X$. \qed
\end{Proof}

\begin{theorem}
Let $\Phi \in \operatorname{Act(G, X)}$ and $H$ be a dense countable subgroup of $G$. Then $P(\Phi)$ is an $F_{\sigma \delta}$-subset of $X$ and hence measurable.
\label{T3.8}
\end{theorem}
\begin{Proof}
For each pair $\epsilon , \delta > 0$, we define $\mathsf{C}(\epsilon , \delta) = \lbrace x\in X \mid \Gamma_{\epsilon}^{x}(\Phi , \Psi) = \phi,$ for some $\Psi \in \operatorname{Act(G, X)}$ satisfying $d_{H}(\Phi, \Psi) \leq \delta \rbrace$. Note that $X\setminus P(\Phi) = \bigcup\limits_{m=1}^{\infty} \bigcap\limits_{n=1}^{\infty}\mathsf{C} (m^{-1} , n^{-1})$. From Lemma \ref{L3.7}, we get that $\mathsf{C}(m^{-1} , n^{-1})$ is an open subset of $X$, for each pair $m , n \in \mathbb{N}$ implying that $X\setminus P(\Phi)$ is a $G_{\delta \sigma}$-subset of $X$ and hence $P(\Phi)$ is an $F_{\sigma \delta}$-subset of $X$. \qed
\end{Proof}

\begin{Definition}
Let $\Phi \in \operatorname{Act(G, X)}$ and $H$ be a dense countable subgroup of $G$. We say that a measure $\mu \in \mathcal{M}(X)$ is a $\beta$-persistent measure (with respect to $H$ and $\Phi$) if for each $\epsilon > 0$, there exists a $\delta > 0$ such that $\mu(B(\epsilon, \Phi, \Psi))  = 1$, for each $\Psi\in \operatorname{Act(G, X)}$ satisfying $d_{H}(\Phi, \Psi) \leq \delta$. The set of all $\beta$-persistent measures (with respect to $H$ and $\Phi$) is denoted by $\mathsf{M}_{P}(\Phi)$.
\label{D3.9}
\end{Definition}

\begin{Remark}
Let $\Phi \in \operatorname{Act(G, X)}$, $H$ be a dense countable subgroup of $G$ and $\mu \in \mathcal{M}(X)$. If for each $\epsilon > 0$, there exists a $\delta > 0$ and a Borelian $B \subseteq X$ with $\mu(B) = 1$ such that for each $\Psi\in \operatorname{Act(G, X)}$ satisfying $d_{H}(\Phi, \Psi) \leq \delta$ and for each $x\in B$, there exists a $y\in X$ such that $d(\Phi_{g}(x), \Psi_{g}(y))\leq \epsilon$, for each $g\in G$, then $\mu\in \mathsf{M}_{P}(\Phi)$.
\label{R3.10}
\end{Remark}

\begin{Definition}
Let $\Phi \in \operatorname{Act(G, X)}$ and $H$ be a dense countable subgroup of $G$. We say that a measure $\mu \in \mathcal{M}(X)$ is an almost $\beta$-persistent measure (with respect to $H$ and $\Phi$) if $\mu(P(\Phi)) = 1$. The set of all almost $\beta$-persistent measures (with respect to $H$ and $\Phi$) is denoted by $\mathsf{M}_{AP}(\Phi)$.
\label{D3.11}
\end{Definition}

\begin{Remark}\label{E1R3.14.1}
Let $\Phi\in \operatorname{Act(G, X)}$, $H_{1}, H_{2}$ be dense countable subgroups of $G$ and $\mu\in \mathcal{M}(X)$. From Remark \ref{R2.1}, we get that $\mu$ is a $\beta$-persistent measure (almost $\beta$-persistent measure, respectively) with respect to $H_{1}$ and $\Phi$ if and only if $\mu$ is a $\beta$-persistent measure (almost $\beta$-persistent measure, respectively) with respect to $H_{2}$ and $\Phi$. Therefore, the notions of $\beta$-persistent measures and almost $\beta$-persistent measures are independent of the choice of dense countable subgroups. 
\end{Remark}

\begin{theorem}
Let $\Phi \in \operatorname{Act(G, X)}$ and $H$ be a dense countable subgroup of $G$. Then the following statements are true:
\begin{enumerate}
\item[(1)] $P(\Phi) = \lbrace x\in X \mid m_{x}\in \mathsf{M}_{P}(\Phi)\rbrace$.
\item[(2)] If $\lbrace \mu_{1},$ . . ., $\mu_{k}\rbrace \subseteq \mathsf{M}_{P}(\Phi)$ and $t_{1},$ . . ., $t_{k}\in (0,1]$ satisfy $\sum_{i=1}^{k}t_{i} = 1$, then $\sum_{i=1}^{k}t_{i}\mu_{i}\in \mathsf{M}_{P}(\Phi)$.
\end{enumerate}
\label{T3.12}
\end{theorem}
\begin{Proof}
We proceed as follows:
\begin{enumerate}
\item[\textit{(1)}] Note that if $x \in P(\Phi)$, then $supp(m_x) = \lbrace x\rbrace \subseteq P(\Phi)$ implying that $m_x \in \mathsf{M}_{P}(\Phi)$. 
Conversely, choose an $m_x \in \mathsf{M}_{P}(\Phi)$, for some $x\in X$. For a given $\epsilon > 0$, choose a $\delta > 0$ such that $m_{x}(B(\epsilon, \Phi, \Psi)) = 1$ whenever $d_{H}(\Phi, \Psi) \leq \delta$ by $\beta$-persistent property of $m_x$ (with respect to $\Phi$). Since $x\in B(\epsilon, \Phi, \Psi)$ and $\epsilon$ is chosen arbitrarily, we get that $x\in P(\Phi)$.
\item[\textit{(2)}] Set $\mu = \sum_{i=1}^{k}t_{i}\mu_{i}$. For a given $\epsilon > 0$, choose $\delta_{1},$ . . . $, \delta_{k} > 0$ by the definition of $\beta$-persistent measures $\mu_{1},$ . . . $, \mu_{k}$ (with respect to $\Phi$) respectively. Set $\delta = \min \lbrace \delta_{1},$ . . . $, \delta_{k}\rbrace$. Note that $\mu(B(\epsilon, \Phi, \Psi)) = \sum_{i=1}^{k} t_{i}\mu_{i}(B(\epsilon, \Phi, \Psi)) = \sum_{i=1}^{k}t_{i} = 1$, for each $\Psi\in Act(G, X)$ satisfying $d_{H}(\Phi, \Psi)\leq \delta$. Since $\epsilon$ is chosen arbitrarily, we get that $\mu \in \mathsf{M}_{P}(\Phi)$. \qed
\end{enumerate}
\end{Proof}

\begin{Lemma}
Let $\Phi\in \operatorname{Act(G, X)}$ and $H$ be a dense countable subgroup of $G$. Then $\Phi$ is $\beta$-persistent if and only if $\mathsf{M}_{P}(\Phi) = \mathcal{M}(X)$.
\label{L3.13}
\end{Lemma}
\begin{Proof}
Note that if $\Phi$ is $\beta$-persistent, then for each $\mu \in \mathcal{M}(X)$ and for each $\epsilon > 0$, there exists a $\delta > 0$ and a Borelian $B = X$ with $\mu(B) = 1$ such that for each $\Psi\in \operatorname{Act(G, X)}$ satisfying $d_{H}(\Phi, \Psi) \leq \delta$ and for each $x\in B$, there exists a $y\in X$ such that $d(\Phi_{g}(x), \Psi_{g}(y))\leq \epsilon$, for each $g\in G$. From Remark \ref{R3.10}, we get that $\mu \in \mathsf{M}_{P}(\Phi)$, for each $\mu \in \mathcal{M}(X)$. 
Conversely, suppose that $\mathsf{M}_{P}(\Phi) = \mathcal{M}(X)$ but $\Phi$ is not $\beta$-persistent. Therefore there exists an $\epsilon > 0$, a sequence of actions $\lbrace \Psi_{i}\rbrace_{i\in \mathbb{N}}\subseteq \operatorname{Act(G, X)}$ and a sequence of elements $\lbrace x_{i}\rbrace \subseteq X$ such that $d_{H}(\Phi , \Psi_{i})\leq \frac{1}{i}$, for each $i\in \mathbb{N}$ and $\Gamma_{\epsilon}^{x_{i}}(\Phi , \Psi_{i}) = \phi$, for each $i\in \mathbb{N}$. Define $\mu = \sum_{i=1}^{\infty}\frac{m_{x_{i}}}{2^{i}}$. For the above $\epsilon$, choose a $\delta > 0$ by $\beta$-persistent property of $\mu$. 
Choose a $k\in \mathbb{N}$ such that $\frac{1}{k} < \delta$. Then there exists a Borelian $B_{k} = B(\epsilon, \Phi, \Psi_{k})$ with $\mu(B_{k}) = 1$ such that $\Gamma_{\epsilon}^{y}(\Phi , \Psi_{k})\neq \phi$, for each $y\in B_{k}$. Since $\mu(B_{k}) = 1$, we get that $m_{x_{k}}(B_{k}) = 1$ implying that $x_{k}\in B_{k}$. Therefore we get that $\Gamma_{\epsilon}^{x_{k}}(\Phi , \Psi_{k}) \neq \phi$, which is a contradiction. \qed
\end{Proof}

\begin{theorem} \label{T3.14}
Let $\Phi \in \operatorname{Act(G, X)}$ and $H$ be a dense countable subgroup of $G$. Then $\mathsf{M}_{P}(\Phi)$ is an $F_{\sigma \delta}$-subset of $\mathcal{M}(X)$.
\end{theorem}
\begin{Proof}
For a given $\epsilon > 0$ and a $\delta > 0$, we denote $\mathsf{C}(\epsilon , \delta) = \lbrace \mu \in \mathcal{M}(X) \mid \mu(B(\epsilon , \Phi , \Psi))$ $< 1, \text{ for some } \Psi\in \operatorname{Act(G, X)} \text{ satisfying } d_{H}(\Phi,\Psi) \leq \delta \rbrace$. Note that $\mathsf{C}(\epsilon , \delta)$ is an open subset of $\mathcal{M}(X)$, for each pair $\epsilon, \delta > 0$ and $\mathcal{M}(X)\setminus \mathsf{M}_{P}(\Phi) = \bigcup\limits_{m=1}^{\infty}\bigcap\limits_{n=1}^{\infty}\mathsf{C}(m^{-1}, n^{-1})$. Therefore $\mathcal{M}(X)\setminus \mathsf{M}_{P}(\Phi)$ is a $G_{\delta \sigma}$-subset of $\mathcal{M}(X)$ and hence $\mathsf{M}_{P}(\Phi)$ is an $F_{\sigma \delta}$- subset of $\mathcal{M}(X)$. \qed
\end{Proof}

\begin{theorem}
Let $\Phi\in \operatorname{Act(G, X)}$ and $H$ be a dense countable subgroup of $G$. If $X$ has no isolated points, then $\Phi$ is $\beta$-persistent if and only if every non-atomic Borel probability measure of $X$ is a $\beta$-persistent measure (with respect to $\Phi$).
\label{T3.15}
\end{theorem}
\begin{Proof}
Forward implication follows from Lemma \ref{L3.13}. Conversely, suppose that every non-atomic Borel probability measure of $X$ is a $\beta$-persistent measure but $\Phi$ is not $\beta$-persistent. Therefore there exists an $\epsilon > 0$, a sequence of actions $\lbrace \Psi_{i}\rbrace_{i\in \mathbb{N}}\subseteq \operatorname{Act(G, X)}$ and a sequence of elements $\lbrace x_{i}\rbrace \subseteq X$ such that $d_{H}(\Phi , \Psi_{i})\leq \frac{1}{i}$, for each $i\in \mathbb{N}$ and $\Gamma_{\epsilon}^{x_{i}}(\Phi , \Psi_{i}) = \phi$, for each $i\in \mathbb{N}$. Recall that if $X$ is a compact metric space without isolated points, then for each $x \in X$ and for each $\epsilon > 0$, there exists a non-atomic Borel probability measure $\nu$ such that $x\in supp(\nu)\subseteq B[x, \epsilon]$ \cite[Lemma 7]{DLMP}. Therefore there exists a sequence of non-atomic Borel probability measures $\lbrace \mu_{i}\rbrace$ such that $x_{i} \in supp(\mu_{i})\subseteq B[x_{i}, \epsilon]$, for each $i \in \mathbb{N}$. Define $\mu = \sum_{i=1}^{\infty}\frac{\mu_{i}}{2^{i}}$. Now, we can complete the proof by using similar arguments as given in the proof of Lemma \ref{L3.13}. \qed
\end{Proof}

\begin{Lemma}
Let $\Phi\in \operatorname{Act(G,X)}$ and $H$ be a dense countable subgroup of $G$. If $\Phi$ is equicontinuous, then $P(\Phi)$ is a closed subset of $X$.
\label{L3.16}
\end{Lemma}
\begin{Proof}
Let $\lbrace x_{i}\rbrace_{i\in \mathbb{N}}$ be a sequence of $\beta$-persistent points of $\Phi$ converging to the point, say $x\in X$. Choose an $\epsilon  > 0$. For $\frac{\epsilon}{2}$, choose an $\eta > 0$ by the equicontinuity of $\Phi$. Choose a $k\in \mathbb{N}$ such that $d(x, x_{k}) < \eta$. Therefore $d(\Phi_{g}(x), \Phi_{g}(x_{k})) < \frac{\epsilon}{2}$, for each $g\in G$. For $\frac{\epsilon}{2}$, choose a $\delta > 0$ by the definition of $\beta$-persistent point $x_{k}$ of $\Phi$. Therefore for each $\Psi\in \operatorname{Act(G,X)}$ satisfying $d_{H}(\Phi, \Psi) \leq \delta$, there exists a $y\in X$ such that $d(\Phi_{g}(x_{k}), \Psi_{g}(y)) \leq \frac{\epsilon}{2}$ implying that $d(\Phi_{g}(x), \Psi_{g}(y)) \leq \epsilon$. Since $\epsilon$ is chosen arbitrarily, we get that $x$ is a $\beta$-persistent point of $\Phi$. \qed
\end{Proof}

\begin{Lemma}
Let $\Phi\in \operatorname{Act(G, X)}$ and $H$ be a dense countable subgroup of $G$. If $\Phi$ is equicontinuous, then $\mathsf{M}_{APe}(\Phi) \subseteq \mathsf{M}_{P}(\Phi)$.
\label{L3.17}
\end{Lemma}
\begin{Proof}
Suppose that $\mu \in \mathsf{M}_{APe}(\Phi) \setminus \mathsf{M}_{P}(\Phi)$. Then there exists an $\epsilon > 0$, a sequence of actions $\lbrace \Psi^{i}\rbrace_{i\in \mathbb{N}} \subseteq \operatorname{Act(G,X)}$ satisfying $d_{H}(\Phi, \Psi^{i})\leq \frac{1}{i}$, for each $i\in \mathbb{N}$ and a sequence of positive measurable sets $\lbrace B_{i}\rbrace_{i\in \mathbb{N}} \subseteq \mathcal{B}(X)$ such that $\Gamma_{\epsilon}^{z}(\Phi , \Psi^{i}) = \phi$, for each $z\in B_{i}$ and for each $i\in \mathbb{N}$. 
Therefore $(B_{i}\cap supp(\mu))\neq \phi$, for each $i\in \mathbb{N}$. Since $\mu\in \mathsf{M}_{APe}(\Phi)$ and $\Phi$ is equicontinuous, we use Lemma \ref{L3.16} to get that $supp(\mu)\subseteq \overline{P(\Phi)} = P(\Phi)$ implying that $(B_{i}\cap supp(\mu))\subseteq (B_{i}\cap P(\Phi))\neq \phi$, for each $i\in \mathbb{N}$. Therefore by compactness of $X$ and from Lemma \ref{L3.16}, we can choose a sequence $\lbrace x_{i} \in B_{i}\cap P(\Phi)\rbrace_{i\in \mathbb{N}}$ converging to some point $x\in P(\Phi)$. 
For $\frac{\epsilon}{2}$, choose a $\delta > 0$, by the equicontinuity of $\Phi$ and $\beta$-persistent property of $\Phi$ through $\lbrace x\rbrace$. Choose a $k\in \mathbb{N}$ such that $\max \lbrace d(x, x_{k}), d_{H}(\Phi, \Psi^{k})\rbrace < \delta$. Therefore $d(\Phi_{g}(x), \Phi_{g}(x_{k})) < \frac{\epsilon}{2}$, for each $g\in G$ and $d(\Phi_{g}(x), \Psi_{g}^{k}(y)) \leq \frac{\epsilon}{2}$, for each $g\in G$ and for some $y\in X$. Hence $d(\Phi_{g}(x_{k}), \Psi^{k}_{g}(y)) \leq \epsilon$, for each $g\in G$ implying that $y\in \Gamma_{\epsilon}^{x_{k}}(\Phi, \Psi^{k})$. Since $x_{k}\in B_{k}$, we get a contradiction and hence $\mu \in \mathsf{M}_{P}(\Phi)$. \qed
\end{Proof}

\textbf{Proof of Theorem \ref{T1.1}.}
Let $H$ be a dense countable subgroup of $G$. First we claim that every equicontinuous topologically stable point of $\Phi$ is a $\beta$-persistent point of $\Phi$. Choose an $x\in E(\Phi) \cap T(\Phi)$ and an $\epsilon > 0$. 
For $\frac{\epsilon}{2}$, choose an $\eta > 0$ by the equicontinuity of $\Phi$ at $x$. For $\eta$, choose a $\delta > 0$ by the definition of topologically stable point $x$ of $\Phi$. Choose a $\Psi \in \operatorname{Act(G, X)}$ satisfying $d_{H}(\Phi, \Psi) \leq \delta$. 
Therefore there exists a continuous map $h: \overline{\mathcal{O}_{\Psi}(x)} \rightarrow X$ such that $\Phi_{g} h = h \Psi_{g}$, for each $g\in G$ and $d(h(z), z) \leq \delta$, for each $z\in \overline{\mathcal{O}_{\Psi}(x)}$. Therefore $d(h(x), x) \leq \delta$ and hence $d(\Phi_{g}(h(x)), \Phi_{g}(x)) \leq \frac{\epsilon}{2}$, for each $g\in G$. 
Hence $d(\Phi_{g}(x), \Psi_{g}(x)) \leq d(\Phi_{g}(x), \Phi_{g}(h(x))) + d(h(\Psi_{g}(x)), \Psi_{g}(x)) \leq \epsilon$, for each $g\in G$. Since $\epsilon$ is chosen arbitrarily, we get that $x\in P(\Phi)$. Since $\Phi$ is an equicontinuous pointwise topologically stable action, we get that $E(\Phi) = T(\Phi) = X = P(\Phi)$ implying that $\Phi$ is pointwise $\beta$-persistent and hence $\mathsf{M}_{APe}(\Phi) = \mathcal{M}(X)$. From Lemma \ref{L3.13} and Lemma \ref{L3.17}, we get that $\Phi$ is $\beta$-persistent. \qed

\begin{Example}
Let $g$ be an equicontinuous pointwise topologically stable homeomorphism of an uncountable perfect compact metric space $(Y,d_0)$ \cite{KTS}. Let $p$ be a periodic point of $g$ with prime period $t\geq 1$. Let $X=Y\cup E$, where $E$ is an infinite enumerable set. Set $Q=\bigcup_{k\in\mathbb{N}} \lbrace 1,2,3\rbrace\times\lbrace k\rbrace\times\lbrace 0,1,2,3,...,t-1\rbrace$. Suppose that $r:\mathbb{N}\rightarrow E$ and $s:Q\rightarrow \mathbb{N}$ are bijections. Consider the bijection $q:Q\rightarrow E$ defined as $q(i,k,j)=r(s(i,k,j))$, for each $(i,k,j)\in Q$. Therefore any point $x\in E$ has the form $x=q(i,k,j)$, for some $(i,k,j)\in Q$. 
%\par 
%\vspace*{0.1cm}

Consider the function $d:X\times X\rightarrow\mathbb{R}^+$ defined by  

\[d(a,b)=\begin{cases} 
0 & \textnormal {if $a=b$},
\\
d_0(a,b) & \textnormal {if $a,b\in Y$} 
\\
\frac{1}{k}+d_0(g^j(p),b) & \textnormal {if $a=q(i,k,j)$ and $b\in Y$}
\\
\frac{1}{k}+d_0(a,g^j(p)) & \textnormal {if $a\in Y$ and $b=q(i,k,j)$} 
\\
\frac{1}{k} & \textnormal {if $a=q(i,k,j)$,$b=q(l,k,j)$ and $i\neq l$} 
\\
\frac{1}{k}+\frac{1}{m}+d_0(g^j(p),g^r(p)) & \textnormal {if $a=q(i,k,j)$,$b=q(i,m,r)$ and $k\neq m$ or $j\neq r$}    
\end{cases}\] 
and $f:X\rightarrow X$ defined by 
\begin{center}
\[f(x)=\begin{cases}
g(x) & \textnormal {if $x\in Y$}
\\
q(i,k,(j+1))$ mod $t & \textnormal {if $x=q(i,k,j).$} 
\end{cases}\] 
\end{center}

Following steps as in \cite[(4) on Page 3742-3743]{CCN}, we get that $(X,d)$ is the compact metric space and $f$ is the homeomorphism.
Note that $x\in X$ is an isolated point if and only if $x\in E$. Let $G = \mathbb{Z}$ be the group of integers. Now, define the action $\Phi\in \operatorname{Act(G, X)}$ by $\Phi_{1}(x) = f(x)$, for each $x\in X$. Since $g$ is equicontinuous, we get that $f$ is equicontinuous and hence $\Phi$ is equicontinuous. 

Choose a $z = q(i,k,j)\in E$ and an $\epsilon < \frac{1}{k}$.  Let $\delta = \epsilon$ and $\Psi \in \operatorname{Act(G, X)}$ be such that $d_{G}(\Phi, \Psi) \leq \delta$. 
Since  $B_{d}(f^{i}(z),\epsilon) = \lbrace f^{i}(z)\rbrace$, for all $-t\leq i\leq t$, we get that $\Phi_{g}(z) = \Psi_{g}(z)$, for each $g\in G$ and $\mathcal{O}_{\Psi}(z) = \lbrace z, f(z),$ . . ., $f^{t-1}(z)\rbrace = \overline{\mathcal{O}_{\Psi}(z)}$. Define $h: \mathcal{O}_{\Psi}(z) \rightarrow X$ such that $h(\Psi_{g}(z)) = \Phi_{g}(z)$, for each $g\in G$. Clearly $h$ is a well defined continuous map such that $\Phi_{g} h = h \Psi_{g}$, for each $g\in G$ and $d(h(x), x) < \epsilon$, for each $x\in \overline{\mathcal{O}_{\Psi}(z)}$. Since $\epsilon$ and $z$ are chosen arbitrarily, we get that $E\subseteq T(\Phi)$.
Now, choose a $y\in Y$ and an $\epsilon > 0$. For this $\epsilon$, choose a $\delta > 0$ by the definition of topologically stable point $y$ of $g$. 
Choose a $\Psi\in Act(G, X)$ satisfying $d_{G}(\Phi, \Psi) \leq \delta$. Since $Y$ has no isolated points, we have $\Psi_{1}(Y) = Y$ and hence $\overline{\mathcal{O}_{\Psi}(y)} = \overline{\mathcal{O}_{\Psi_{1}}(y)} \subseteq Y$. Therefore we can choose a continuous map $h: \overline{\mathcal{O}_{\Psi}(y)} \rightarrow X$ such that $\Phi_{g} h = h \Psi_{g}$, for each $g\in G$ and $d(h(x), x) < \epsilon$, for each $x\in$ $\overline{\mathcal{O}_{\Psi}(y)}$. Since $y$ and $\epsilon$ are chosen arbitrarily, we get that $Y\subseteq T(\Phi)$  implying that $\Phi$ is an equicontinuous pointwise topologically stable action. From Theorem \ref{T1.1}, we get that $\Phi$ is $\beta$-persistent.
\label{E3.18}
\end{Example}

\begin{Remark}\label{R3.19}
Let $G$ be a finitely generated group (not necessarily first countable and Hausdorff) and $\Phi\in \operatorname{Act(G, X)}$. Following similar steps as in the proof of Theorem \ref{T1.1}, we can prove that if $\Phi$ is pointwise topologically stable in the sense of \cite{DJLP}, then $\Phi$ is $\beta$-persistent in the sense of \cite{ALLP}. 
\end{Remark}

\textbf{Acknowledgements:} The first author is supported by CSIR-Senior Research Fellowship (File No.-09/045(1558)/2018-EMR-I) of  Government of India.

\end{document}